\documentclass[11pt]{article}
\newtheorem{theorem}{Theorem}
\newtheorem{lemma}{Lemma}

\def\katsu{\mathop{\wedge}}
\def\cov{{\rm cov}}
\def\F{{\cal F}}
\def\G{{\cal G}}
\def\s{{\cal S}}
\def\Y{{\cal Y}}
\def\x{{\bf x}}
\def\b{{\bf b}}
\def\<{\langle}
\def\>{\rangle}

\title{
Semiparametric estimation of volatility functions of diffusion processes
from discretely observed data
}
\author{
Isao Shoji
\\
Graduate School of Systems and Information Engineering,\\
University of Tsukuba,
Tsukuba Ibaraki, 305-8573, Japan\\
{\rm shoji@sk.tsukuba.ac.jp}
}

\begin{document}

\maketitle

\begin{abstract}
This paper provides a semiparametric model to estimate processes of the volatility defined as the squared diffusion coefficient of a stochastic differential equation. Without assuming any functional form of the volatility function, we estimate the volatility process by filtering.
We prove the consistency of the model in the sense that estimated processes converge to the true ones as the number of observations ($N$) goes to infinity and the sampling time interval ($\Delta t$) goes to zero while $N\Delta t$ going to infinity.
We also carry out numerical experiments through stochastic differential equations with linear/nonlinear volatility functions in order to check whether or not the model can actually estimate the volatility and compare the performance with the local linear model.




\end{abstract}


\section{Introduction}
When modeling time series by continuous-time stochastic processes,
we often face a difficult problem of what kind of functions should be used for the drift and diffusion coefficients of a stochastic differential equation since we have almost no knowledge about them beforehand.
But, the specification of the diffusion coefficient is much more important for the modeling. Actually, recent researches about analysis of financial time series show the weak evidence of nonlinearity in the drift suggested by Stanton (1997) for example; Chapman and Peason (2000) addresses that the test of the nonlinearity is not robust through the simulation studies. And, Fan and Zhang (2003) develops an alternative test free from the problem of the method used by Stanton (1997) and show the weak evidence against the linear drift of Standard \& Poor 500 as well as the short-term interest rate. Furthermore, Sun (2003) and recently Bali and Wu (2006) report the similar results.
Additionally, from a technical point of views as pointed out by Bandi and Phillips (2003), the drift coefficient cannot be identified nonparametrically on a fixed time interval.

To the contrary, those researches stress the nonlinearity in the diffusion coefficient, or the volatility, which is crucial for describing the time evolution of financial time series such as interest rate data.
And besides, there is no such a technical problem of identification as estimating the drift coefficient on a fixed interval.
So, the specification of the volatility is really important for their modeling.

In the recent statistical models of volatility, the realized volatility is 
becoming one of the most successful tools in modeling and forecasting the 
volatility, and its more extending studies have been extensively carried out 
recently; Thomakos et al (2002), Deo et al (2006), Engle and Gallo (2006), and Ghysels et al (2006), for example. The realized volatility is basically rooted in the fact that the quadratic variation on a time interval converges in probability to the volatility over that time interval, so-called the integrated volatility. Theoretical and numerical studies such as Andersen et al (2003, 2004, 2005) and Barndorff-Nielsen and Shephard (2002, 2004) that are relevant to stochastic volatility model show that the estimation by the realized volatility is well performed. Though the realized volatility can effectively estimate the integrated volatility, it's still difficult to estimate the spot volatility, or the squared diffusion coefficient, that is defined as the integrand of the integrated volatility. This information is indispensable for setting up a stochastic differential equation and using it for practical purposes as well. And besides, the integrated volatility easily recovers from the spot volatility, but it's not easy to do the converse.
To estimate the spot volatility, we usually need to have some information about its functional form beforehand. But, this is not expected because of little knowledge about the functional form of the spot volatility.

The aim of the paper is to present a method of estimating the spot volatility, simply called the volatility in the paper, of a one-dimensional stochastic differential equation from discrete observations. But, since we have no knowledge about what kind of functions should be used for the volatility, we need to model it nonparametrically. The most straightforward way is to use its polynomial approximation. But this approach doesn't seem successful since estimation of a polynomial function is not efficient particularly when a higher order polynomial is used. Conversely, however, we can't use a lower order polynomial since it leads to bad approximation of the diffusion coefficient after all.

The drawbacks of this approach lie in fitting a polynomial globally. Hence, we could use local polynomial modeling as alternatives. This modeling is based on the kernel regression in which the regression function is expressed as the weighted average of several sub-regression functions that are usually first or second order polynomials and these weighted are characterized by the so-called kernel function. See Fan and Gijbels (1996) and Campbell et al (1997), for example.
Actually, the kernel regression, more specifically the local polynomial model, is used for estimating the volatility in a nonparametric manner from Florens-Zmirou (1993) to Stanton (1997), Fan and Yao (1998), Jacod (2000), Bandi and Phillips (2003) and Fan and Zhang (2003), which are given as fully nonparametric models while A\"{\i}t-Sahalia (1996) proposes a semiparametric model in which the functional form of the drift coefficient is known.
Though the local polynomial modeling doesn't suffer from the trouble of higher order, as pointed out in Campbell et al (1997), it instead has the problem of overfitting and bandwidth selection. In particular the overfitting is serious in forecasting the volatility.

To avoid these intractability, we reconsider the local polynomial modeling from a different point of views. In the local polynomial modeling, though it is considered as a nonparametric model, each polynomial over its window needs to be estimated parametrically. But, this paper proposes a model of estimating processes of the volatility function at observed states of the process without estimating its parametric functional form. Simply stated, every unobservable process of the volatility is constructed out of the observable ones.
Or intuitively, the model looks like a local polynomial model with infinitesimal bandwidth.
And then, we try to get a one-to-one correspondence between the observable and unobservable processes and draw them on the plane, which will produce information on the functional form of the volatility.

This method depends solely on how unobservable processes should be estimated from 
observable ones. The state space modeling is one of the most popular methods for that purpose so that every unobservable state can be easily estimated from observable states thanks to the Kalman filtering. So it seems we have only to set up a state space model in which we define states of the volatility as unobservable ones. 
Although the method is surely straightforward, we can't directly apply the updating formula of the Kalman filtering to the problem under consideration since a stochastic differential equation cannot be necessarily handled by its simple application. In this paper, we propose an alternative recursive updating formula, and thereby we get estimates of the volatility as filtered states.
Thanks to the recursive updating, the prediction and filtering can reflect the recent state of the process so that we don't have to care much about the problem of the overfitting. And besides, no bandwidth selection is required. 

From a theoretical viewpoint, it is quite important whether or not the proposed model has the consistency in the sense that estimated processes converge to the true ones as sampling time interval goes to zero for example. We present a proof of the consistency so that theoretically we can estimate true processes as accurately as we might expect by making the sampling time interval close to zero while making the total time span as large as possible.
On the other hand, from a practical point of views, it is important as well whether we can feasibly implement the model or not.
By using stochastic differential equations with linear/nonlinear volatilities, 
we carry out numerical experiments to see how well we can estimate volatility functions from discretely observed data. And, we compare the performance of the proposed model with the local linear model which is used as one of the local polynomial models.
Additionally, we estimate the integrated volatilities by using the estimated volatility processes and compare them with those estimated by the realized volatility.

The organization of this paper is as follows. Firstly we propose a model by which unobservable processes of the volatility function can be estimated from discrete time series of the process of interest. Secondly, we discuss the consistency of the model by investigating the asymptotic behavior of estimated processes. And then, we conduct Monte Carlo experiments to evaluate the performance of the model through comparison with the local linear model. Last, we give the concluding remarks.

\section{Semiparametric Model}\label{model}
We consider a one-dimensional diffusion process, $X_t$, which never explodes in finite time and satisfies the following stochastic differential equation (SDE) starting at a constant $\xi>0$,
\begin{equation}\label{SDE}
dX_t=\mu(X_t;\eta)dt+\sigma(X_t) dB_t,
\end{equation}
where $\mu(x;\eta)$ is a linear/nonlinear function which is twice continuously differentiable with respect to $x$ and $\eta$ and $\{B_t,\F_t\}_{t\ge 0}$ is a standard Brownian motion on filtration $\{\F_t\}_{t\ge 0}$. We define the volatility function $g(x)$ by $g(x)\equiv \sigma(x)^2$ and assume $\int_0^tg(X_u)du<\infty$ for any $t<\infty$ almost surely.
As in the same setting as A\"{\i}t-Sahalia (1996), we assume $\mu(x;\eta)$ is known with an unknown parameter vector $\eta$ but $g(x)$ is completely unknown.
Instead, we could assume $\mu$ is completely unknown as well if we have a method to estimate consistently $\mu$ and $\mu'$ with its rate of convergence $(N\Delta t)^{-{1\over 2}}$, where $N$ is the number of observations and $\Delta t$ is a sampling interval. But, throughout this paper, we assume $\mu(x;\eta)$ is known while $\eta$ is unknown.

Suppose equidistant discrete times expressed by $0=t_0<t_1<\cdots<t_N=T$ with $\Delta t=T/N$.
We observe the process $X_t$ at the discrete times, $\{X_{t_k}\}_{1\le k\le N}$.
Under this setting, we want to estimate the discrete states
$\{g(X_{t_k})\}_{1\le k\le N}$ from $\{X_{t_k}\}_{1\le k\le N}$.

First, suppose an approximation of $g$, denoted by $f$, given as the second order Taylor expansion around $x_0$:
\begin{equation}\label{def:f}
f(x)\equiv g(x_0)+g'(x_0)(x-x_0)+{g''(x_0)\over 2}(x-x_0)^2
\end{equation}
By replacing $x$ by the process $X_t$, we can approximate 
$g(X_t)$ as a quadratic function of $X_t$. Here assuming $x_0$ to be fixed globally, the approximation leads to a global polynomial approximation. Instead, replacing $x_0$ by $X_s$ which changes depending upon choice of $s$, we get a local polynomial approximation.
In the local polynomial approximation, the coefficients such as $g'$ and $g''$ are constant over $[t_{k-1},t_k)$ where $s=t_{k-1}$, but not globally. Hence, even if $g(x)$ is actually a cubic function for example, it could be well approximated by the local polynomial model of degree two just as a smooth curve can be approximated piecewisely by tangent lines.
By contrast, the global polynomial model frequently comes to bad approximation 
particularly when $g$ shows high nonlinearity.

In the local polynomial approximation, we define new processes, $Y_t$, $Y^1_t$, and $Y^2_t$, by
\begin{eqnarray*}
Y_t&\equiv& f(X_t)\\
Y^1_t&\equiv& f'(X_t)\\
Y^2_t&\equiv& f''(X_t).
\end{eqnarray*}
In order to see how these processes evolve in time, we
apply the It\^o's formula to $Y_t$, $Y^1_t$, and $Y^2_t$ on $t_{k-1}\le s<t<t_k$, and thereby we get,
\begin{eqnarray*}
Y_t-Y_s&=& \int_s^tY^1_udX_u+\int_s^t{1\over 2}Y^2_ud\<X\>_u\\
Y^1_t-Y^1_s&=& \int_s^tY^2_udX_u\\
Y^2_t-Y^2_s&=& 0.
\end{eqnarray*}
The last equality implies $Y^2_t$ is constant over $[t_{k-1},t_k)$. But, we proceed as if $Y^2_t$ to be globally constant
and denote it by $\theta$ in place of $Y^2_t$. Using this, we rewrite the above system in a differential form as follows:
\begin{eqnarray}
dY_t&=&Y^1_tdX_t+{\theta\over 2} d\<X\>_t\label{SDE:Y}\\
dY^1_t&=&\theta dX_t\label{SDE:Y1}
\end{eqnarray}
under the original SDE (\ref{SDE}). Here we set $Y_s=g(X_s)$ at every 
$s\in\{t_k\}_{0\le k\le N}$.
Hence, every sample path of $Y_t$ is not necessarily continuous at
$\{t_k\}_{0\le k\le N}$, but is continuous over $[t_{k-1},t_k)$ for all $k$
$(1\le k\le N)$.

Combining (\ref{SDE}), (\ref{SDE:Y}) and (\ref{SDE:Y1}), we can setup a system of SDEs. On the assumption, we can observe $X_t$ but not $Y_t$. So, we have to estimate $Y_t$ by the technique of filtering or something like that. The system, however, is not so tractable for the purpose as the linear system which can produce an estimate of such an unobservable process as $Y_t$ through the Kalman-Bucy filtering for example. So, we want another system as an approximation of the system (\ref{SDE}), (\ref{SDE:Y}) and (\ref{SDE:Y1}).

First, we consider a linear approximation of $\mu$ in (\ref{SDE}) around $x_0$, denoting by $\tilde\mu(x;\eta)$ as follows:
\begin{eqnarray*}
\tilde\mu(x;\eta)&\equiv&
\mu(x_0;\eta)+\mu'(x_0;\eta)(x-x_0)\\
&=&
\mu(x_0;\eta)-\mu'(x_0;\eta)x_0+\mu'(x_0;\eta)x.
\end{eqnarray*}
Replace $x$, $x_0$ and $\eta$ by $X_t$, $X_s$ and some estimate of $\eta$, $\hat\eta$, respectively; $\hat\eta$ will be replaced later by a least squar estimate. And, denote $\mu(X_s;\hat\eta)-\mu'(X_s;\hat\eta)X_s$ and $\mu'(X_s;\hat\eta)X_t$ by $\alpha_s$ and $\beta_s$ for simplicity.
Similarly, we replace $Y^1_t$ in (\ref{SDE:Y}) by $Y^1_s$ like the Euler method. And then, to link the observable process and the volatility process as an unobservable one, we define $\tilde X_t$ and $\tilde Y_t$, as approximation of $X_t$ and $Y_t$, which satisfy the following system of SDEs:
\begin{eqnarray}
d\tilde X_t&=&(\alpha_s+\beta_s\tilde X_t)dt+\sqrt{\tilde Y_t}dB_t\label{aSDE:X}\\
d\tilde Y_t&=&Y^1_sd\tilde X_t+{\theta\over 2} d\<\tilde X\>_t\label{aSDE:Y}\\
dY^1_t&=&\theta dX_t
\end{eqnarray}
for $t\in[t_{k-1},t_k)$ $(1\le k\le n)$. And, we take $\tilde X_{t_{k-1}}=X_{t_{k-1}}$ at the end-point. That is, like $Y_t$, we reset the initial state of the approximate observable process to that of the original one at discrete times. As for
$\tilde Y_{t_{k-1}}$, it's recursively defined. Initially $\tilde Y_0=Y_0$. And then, we define it by $\tilde Y_{t_{k-1}}=\lim_{s\uparrow t_{k-1}}\tilde Y_s$. Thereby $\tilde Y_t$ is a continuous process.

Here note the difference between $Y_t$ and $\tilde Y_t$ as well as $X_t$ and $\tilde X_t$. And besides, $Y^1_t$ is the same between the two systems. Firstly, $X_t$ and $Y_t$ evolve in time according to (\ref{SDE}) and (\ref{SDE:Y}), whereas $\tilde X_t$ and $\tilde Y_t$ do according to (\ref{aSDE:X}) and (\ref{aSDE:Y}), respectively.
So, $X_t$ is continuous while $\tilde X_t$ is not necessarily continuous. Conversely, $Y_t$ is not necessarily continuous while $\tilde Y_t$ is continuous.
But, $Y^1_t$ is driven by $X_t$ for the both cases. Actually, $Y^1_t$ is immediately given by $Y^1_t=Y^1_s+\theta(X_t-X_s)$, or equivalently,
$Y^1_t=Y^1_0+\theta(X_t-X_0)$.

Though the system (\ref{aSDE:X})-(\ref{aSDE:Y}) looks like a stochastic volatility model, it differs since the system is derived from stochastic differential equations with time-homogeneous drift and diffusion coefficients. Differently from stochastic volatility models,
the system (\ref{aSDE:X})-(\ref{aSDE:Y}) is tractable since the drift coefficients are locally linear in $\tilde X$ and $\tilde Y$. Hence, the system of SDE's can be explicitly solved on $\F_s$, and their conditional expectations with respect to $\F_s$ are easily obtained. To this end, we rewrite the system compactly as follows:
\begin{equation}
d\x_t=(A\x_t+\b)dt+\s(\x_t)dB_t
\end{equation}
where,
\begin{eqnarray*}
\x_t&=&(\tilde X_t,\tilde Y_t)'\\
A&=&
\left(
\begin{array}{cc}
\beta_s& 0\\
\beta_s Y^1_s& \theta/2
\end{array}
\right),\ 
\b\ =\ 
\left(
\begin{array}{c}
\alpha_s\\
\alpha_s Y^1_s
\end{array}
\right),\ 
\s(\x_t)\ =\ 
\left(
\begin{array}{c}
\sqrt{\tilde Y_t}\\
Y^1_s\sqrt{\tilde Y_t}
\end{array}
\right)
\end{eqnarray*}
For $t_{k-1}\le s<t<t_k$ and $\theta\neq 0$, its solution is given as,
\begin{eqnarray}
\x_t
&=&\exp(A(t-s))\x_s
+A^{-1}(\exp(A(t-s))-I)\b\\
&&+\int_s^t\exp(A(t-u))\s(\x_u) dB_u
\nonumber
\end{eqnarray}
where,
\[
\exp(At)=
\left(
\begin{array}{cc}
e^{at}& 0\\
{b\over c-a}e^{ct}
+{b\over a}\left(1-{c\over c-a}\right)e^{at}& e^{ct}
\end{array}
\right)
\]
with $a=\beta_s$, $b=\beta_s Y^1_s$, and $c=\theta/2$ and $I$ is an identity matrix. Hence, the conditional mean and variance on $\F_s$, denoted by $E_s[\x_t]$ and $\cov_s(\x_t)$ respectively, are given as,
\begin{eqnarray}
E_s[\x_t]
&=&\exp(A(t-s))\x_s
+A^{-1}(\exp(A(t-s))-I)\b\label{c-mv}\\
\cov_s(\x_t)
&=&
E\left[\left.
\int_s^t\exp(A(t-u))\s(\x_u)\s(\x_u)'\exp(A(t-u))'du
\right|\F_s\right]\\
&=&
\int_s^t\exp(A(t-u))E[\s(\x_u)\s(\x_u)'|\F_s]\exp(A(t-u))'du
\nonumber
\end{eqnarray}
Here note that the conditional mean of $\x_t$ is linear in $\x_s$ and all the 
components of $A$ and $\b$ are characterized by the local/global constants, $\alpha_s$, 
$\beta_s$, $\theta$, and $Y^1_s$.
Since
\begin{eqnarray*}
E[\s(\x_u)\s(\x_u)'|\F_s]
&=&
E[\tilde Y_u|\F_s]
\left(
\begin{array}{cc}
1& Y^1_s\\
Y^1_s& (Y^1_s)^2
\end{array}
\right)
\end{eqnarray*}
$\cov_s(\x_t)$ can be further computed by using the formula of $E_s[\x_u]$. After somewhat cumbersome computation, we get,
\begin{equation}
\cov_s(\x_t)=
\left(
\begin{array}{cc}
I_1& pI_1+(q+Y^1_s)I_2\\
pI_1+(q+Y^1_s)I_2& p^2I_1+2p(q+Y^1_s)I_2+(q+Y^1_s)^2I_3
\end{array}
\right)
\end{equation}
where,
\begin{eqnarray*}
I_1&=&
\Y_1e^{a\Delta t}{e^{a\Delta t}-1\over a}
+\Y_2e^{c\Delta t}{e^{(2a-c)\Delta t}-1\over 2a-c}
+\Y_3{e^{2a\Delta t}-1\over 2a}\\
I_2&=&
\Y_1e^{a\Delta t}{e^{c\Delta t}-1\over c}
+\Y_2e^{c\Delta t}{e^{a\Delta t}-1\over a}
+\Y_3{e^{(a+c)\Delta t}-1\over a+c}\\
I_3&=&
\Y_1e^{a\Delta t}{e^{(2c-a)\Delta t}-1\over 2c-a}
+\Y_2e^{c\Delta t}{e^{c\Delta t}-1\over c}
+\Y_3{e^{2c\Delta t}-1\over 2c}\\
\Y_1&=&
p\tilde X_s-{b\alpha_s\over ac}\\
\Y_2&=&
q\tilde X_s+\tilde Y_s+{\alpha_s Y^1_s\over c}\\
\Y_3&=&
{b\alpha_s\over ac}-{\alpha_s Y^1_s\over c}\\
p&=&{b\over a}\left(1-{c\over c-a}\right),
\ q\ =\ {b\over c-a}
\end{eqnarray*}
Since $\Y_2$ is linear in $\tilde Y_s$, $\cov_s(\x_t)$ is a linear function of $\tilde Y_s$ and denoted by $V_{t|s}(\tilde Y_s)$ for simplicity.

These expectations are not immediately used for the estimate of $\tilde Y_t$ that is considered as an approximation of $g_t$, defined by $g_t\equiv g(X_t)$, since the expectations are conditional on $\F_s$. But, we want to estimate it on the condition of the discrete observations.
Let $\G_{t_k}$ be a $\sigma$-algebra generated by $\{X_{t_j}\}_{0\le j\le k}$, and abbreviate $t$ and $s$ for $t_k$ and $t_{k-1}$, respectively.
To this end, we construct estimators of $\x_t$ and $\x_s$ on the condition of $\G_s$, denoted by $\x_{t|s}$ and $\x_{s|s}$, as follows:
\begin{eqnarray}
\x_{t|s}&=&
\exp(A(t-s))\x_{s|s}+A^{-1}(\exp(A(t-s))-I)\b\label{pred:Y}\\
\x_{t|s}&=&(\tilde X_{t|s},\tilde Y_{t|s})'\\
\x_{s|s}&=&(\tilde X_{s|s},\tilde Y_{s|s})'
\end{eqnarray}
Here note $\tilde X_{s|s}$ belongs to $\G_s$ since $\tilde X_s=X_s$ by the setting. 
Hence $\tilde X_{t|s}=E[\tilde X_t|\G_s]=E[\tilde X_t|\F_s]$ since $A$ and $\b$ belong to $\G_s\subset\F_s$.
For $\tilde Y_{t|t}$, we construct an estimator of $\tilde Y_t$ on the condition of $\G_t$ by,
\begin{eqnarray}
\tilde Y_{t|t}&=&
\tilde Y_{t|s}+\kappa(\tilde X_t-\tilde X_s)\label{filter:Y}\\
\kappa&=&
{V_2(\tilde Y_{s|s})
\over V_1(\tilde Y_{s|s})}
\end{eqnarray}
where $V_1$ and $V_2$ are the (1,1) and (1,2) elements of $V_{t|s}$.
By the formula of (\ref{pred:Y}) and (\ref{filter:Y}), $\tilde Y_{t|s}$ and $\tilde Y_{t|t}$ are recursively updated after the initial state is given by $\tilde Y_{0|0}=Y_0$.
Thanks to the recursive formula, $\tilde Y_{t|t}\in\G_t$ for all $t\in\{t_k\}_{0\le k\le N}$. Actually, $\tilde Y_{0|0}$ is known. Suppose $\tilde Y_{s|s}\in\G_s$. Then, $\tilde Y_{t|s}\in\G_s$ due to (\ref{pred:Y}). But, by (\ref{filter:Y}), $\tilde Y_{t|t}\in\G_t$.

Here note $V_{t|s}(\tilde Y_{s|s})$ which isn't necessarily equal to $V_{t|s}(\tilde Y_s)$ belongs to $\G_s$ since all the associated coefficients belong to $\G_s$. And, these formula can be regarded as the prediction and filtering in the Kalman filtering if the system (\ref{aSDE:X})-(\ref{aSDE:Y}) is a conventional linear system.

For practical purpose, we need to know the parameter vector, $\eta$, and the nuisance parameter, $\theta$. $\eta$ can be estimated by the least square estimation for example.
As for $\theta$, we can take any nonzero number regardless of which the consistency discussed in the next section still holds thanks to theorem $\ref{th:gY}$. But, for numerical efficiency, we can take a quasi-maximum likelihood estimate obtained from maximizing the following likelihood function:
\begin{eqnarray}\label{logl}
p(X_{t_0},X_{t_1},\cdots,X_{t_n})&=&
p(X_{t_0})
\prod_{k=1}^n
(
2\pi H
V_{t_k|t_{k-1}}(\tilde Y_{t_{k-1}|t_{k-1}})
H')^{-1/2}\\
\nonumber
&&\times
\exp\left\{-{(X_{t_k}-H\x_{t_k|t_{k-1}})^2
\over 2H
V_{t_k|t_{k-1}}(\tilde Y_{t_{k-1}|t_{k-1}})
H'}\right\}\\
H&=&(1,0)
\nonumber
\end{eqnarray}

\section{Consistency of the model}
In the first place, we preliminarily set up the followings:

\begin{enumerate}
\item {\em Observation}:
Suppose equidistant discrete times expressed by $0=t_0<t_1<\cdots<t_N=T$ with $\Delta t=T/N$.
Let time $\tau$ be arbitrarily fixed.
But, the discrete times necessarily traverse $\tau$; that is, $t_n=\tau$ for some $n$.
The process $X_t$ is observed at the discrete times and denote the observations by
$\{X_{t_k}\}_{0\le k\le N}$.
\item {\em Lipschitz condition}:
$\mu(x;\eta)$ is twice continuously differentiable with respect to $x$ and $\eta$, and
$g(x)$ and $\sqrt{g(x)}$ are twice continuously differentiable as well.
$\mu$ and $\sqrt{g}$ satisfy the Lipschitz conditions. That is, there is a constant $L$ such that,
\begin{eqnarray}
|\mu(x)-\mu(y)|&\le& L|x-y|\\
|\sqrt{g(x)}-\sqrt{g(y)}|&\le& L|x-y|
\end{eqnarray}
\item {\em Localization}:
First we assume $X_t$ and $\tilde Y_t$ never explode in finite time and $\tilde Y_t$ never reaches zero as well.
Suppose a positive number $M$ which is arbitrarily given.
By using the stopping time $T=\katsu_{1\le i\le 3}T_i$ for $T_i$ given below, we define stopped processes such as $X_t\equiv X_{t\katsu T}$,
$\tilde X_t\equiv \tilde X_{t\katsu T}$,
$Y_t\equiv Y_{t\katsu T}$, and
$\tilde Y_t\equiv \tilde Y_{t\katsu T}$, where
\begin{eqnarray}
T_1&=& \inf\{t\ge 0; X_t\ge M\ {\rm or}\ \<X\>_t\ge M\},\\
T_2&=& \inf\{t\ge 0; \tilde Y_t\ge M\},\\
T_3&=& \inf\{t\ge 0; \tilde Y_t\le 0\}.
\end{eqnarray}
Note $X_t$ and $\tilde Y_t$ are continuous processes, and thereby $T_i$'s 
are all suitably defined as stopping time; see Karatzas and Shreve (1991), for example.
By the above definition, we can assume $X_t$ and $\tilde Y_t$ are bounded. And 
besides, $Y_t$ and $Y^1_t$ can be assumed to be bounded, too.
Actually, from (\ref{def:f}),
$Y_t=Y_s+Y^1_s(X_t-X_s)+(\theta/2)(X_t-X_s)^2$.
But, $Y^1_t=Y^1_s+\theta(X_t-X_s)$, or $Y^1_t=Y^1_0+\theta(X_t-X_0)$.
Hence, $Y^1_t$ is bounded, and so is $Y_t$.
Here note $Y_s=g(X_s)$.

By the localization, we first assume $X_t$, $Y_t$, and $\tilde Y_t$ are all bounded, and thereby we prove the following theorems for the bounded processes. And then, by letting $M\to\infty$, we get the final result.
\item {\em Initial state}: The initial states $X_0$, $Y_0$, $Y^1_0$, and $\tilde 
Y_0$ are given as constant. Particularly, we assume $Y_0=\tilde Y_0=\tilde Y_{0|0}$.
\item {\em Asymptotics}: We consider $N\to\infty$, $\Delta t\to0$ and $N\Delta t\to\infty$ simultaneously. Here, $\Delta t=T/N$.
\item {\em Consistent estimate $\hat\eta$}: We assume we can consistently estimate $\eta$ with its rate of convergence $(N\Delta t)^{-{1\over 2}}$. We can take the least square estimation as such a method for example.
See Prakasa Rao (1983).
Another estimators of drift coefficients are known to have the same rate of convergence; see Florens-Zmirou (1989), Yoshida (1992) and Kessler (1997).
\end{enumerate}
Under the above conditions, we want to show the model has the consistency in the sense that $E|g_t-\tilde Y_{t|t}|^2$ converges to zero as $N\to\infty$, $\Delta t\to 0$ and $N\Delta t\to\infty$.
Instead of evaluating directly the measure, we do separately the distance between $g_t$ and $Y_t$ and between $Y_t$ and $\tilde Y_t$. Here recall $Y_t$ evolve in time as follows:
\begin{eqnarray*}
dX_t&=&\mu(X_t;\eta)dt+\sqrt{g_t}dB_t\\
dY_t&=&Y^1_tdX_t+{\theta\over 2}d\<X\>_t\\
dY^1_t&=&\theta dX_t
\end{eqnarray*}
On the other hand, $\tilde Y_t$ follows the system:
\begin{eqnarray*}
d\tilde X_t&=& \tilde\mu(\tilde X_t;\hat\eta)dt+\sqrt{\tilde Y_t}dB_t\\
d\tilde Y_t&=& Y^1_sd\tilde X_t+{\theta\over 2}d\<\tilde X\>_t\\
dY^1_t&=&\theta dX_t
\end{eqnarray*}
Since
\begin{eqnarray*}
(g_t-\tilde Y_{t|t})^2&=&
\{(g_t-Y_t)+(Y_t-\tilde Y_t)+(\tilde Y_t-\tilde Y_{t|t})\}^2,\\
&\le&
3\{(g_t-Y_t)^2+(Y_t-\tilde Y_t)^2+(\tilde Y_t-\tilde Y_{t|t})^2\},
\end{eqnarray*}
we have only to show $E|g_t-Y_t|^2$, $E|Y_t-\tilde Y_t|^2$, and
$E|\tilde Y_t-\tilde Y_{t|t}|^2$ converge to zero.
Actually, we can show the following theorems.
\begin{theorem}\label{th:gY}
For any $t\in\{t_k\}_{0\le k\le n}$, $\lim_{\Delta t\to 0}E|g_t-Y_t|^2=0$.
\end{theorem}
\noindent
{\bf proof:} Thanks to lemma \ref{lm:EXm} and \ref{lm:gY} in the appendix, we can immediately show it.

\noindent
Theorem \ref{th:gY} implies hat $Y_t$ converges to $g_t$ independent of $\theta$. So, theoretically, we don't have to care about its asymptotic properties as far as the consistency of the proposed model is concerned.
\begin{theorem}\label{th:YaY}
Let $\tau$ be arbitrarily fixed. Suppose equidistant discrete times traversing $\tau$; that is, $0=t_0<t_1<\cdots<t_n=\tau<\cdots<t_N=T$ for some $n$.
Then, $E|Y_\tau-\tilde Y_\tau|^2\to0$
as $N\to\infty$, $\Delta t\to 0$ and $N\Delta t\to\infty$.
\end{theorem}
\noindent
{\bf proof:} See the appendix.
\begin{theorem}\label{th:aYfY}
In the same setting as theorem \ref{th:YaY},
$E|\tilde Y_\tau-\tilde Y_{\tau|\tau}|^2\to0$
as $N\to\infty$, $\Delta t\to 0$ and $N\Delta t\to\infty$.
\end{theorem}
\noindent
{\bf proof:} See the appendix.

\medskip
\noindent
Then, we finally get,
\begin{theorem}
In the same setting as theorem \ref{th:YaY},
$E|g_\tau-\tilde Y_{\tau|\tau}|^2\to0$
as $N\to\infty$, $\Delta t\to 0$ and $N\Delta t\to\infty$.
\end{theorem}

\section{Numerical experiments}
First, we try to estimate curves of volatility functions by plotting tuples of $X_t$ and $\tilde Y_{t|t}$ that are estimated by the proposed model. Section 3 guarantees the consistency of the proposed model, so we want to confirm this numerically by seeing how the estimates behave as the sampling interval goes to zero.

Next, we compare the performance of the proposed model with the local linear model, or the local polynomial model when a linear function being fitted locally. According to Fan and Zhang (2003), the local linear model used here is briefly explained as follows.

Taking $m(x)$ for the volatility function and considering the neighborhood of $x_0$, $m(x)$ is locally approximated by $m(x)\approx\hat m(x)\equiv \beta_0+\beta_1(x-x_0)$, where the coefficients $\beta_0$ and $\beta_1$ are given by minimizing the object function,
\[
\sum_{k=1}^n\{Z^*_{t_k}-\beta_0-\beta_1(X_{t_{k-1}}-x_0)\}^2K_h(X_{t_{k-1}}-x_0).
\]
Here $Z^*_{t_k}=(X_{t_k}-X_{t_{k-1}})^2/\Delta t$ and $K_h(\cdot)=K(\cdot/h)/h$. where $\Delta t$ is the sampling interval,  $K(\cdot)$ is a kernel function and $h$ is a bandwidth. We use the Epanechnikov kernel defined by $K(u)=(3/4)(1-u^2)I(|u|\le 1)$ where $I(\cdot)$ is the indicator function.

Throughout the numerical experiments, we consider the case in which $\mu$ is linear since the least square estimation (LSE) produces the consistent estimate of $\eta$ with its rate of convergence $(N\Delta t)^{-{1\over 2}}$. Let $\mu(x)=\alpha+\beta x$ where $\eta=(\alpha,\beta)$. Then, $E[X_t|\F_t]=X_s+(\alpha/\beta+X_s)(\exp(\beta(t-s))-1)$. So, we can get the estimates of $\alpha$ and $\beta$ by minimizing
$\sum_{k=1}^n(X_t-E[X_t|\F_t])^2$ with respect to $\alpha$ and $\beta$. On the other hand, we estimate the nuisance parameter $\theta$ by using quasi-MLE with (\ref{logl}).

\subsection{Estimation of volatility process}

Here, we consider the following SDEs:
\begin{eqnarray*}
dX_t&=&(1-X_t)dt+\sqrt{X_t}dB_t\\
dX_t&=&(1-X_t)dt+\sqrt{X_t^2}dB_t\\
dX_t&=&(1-X_t)dt+\sqrt{X_t^3}dB_t\\
dX_t&=&(1-X_t)dt+\sqrt{X_t\exp(-X_t^2)}dB_t
\end{eqnarray*}
with $X$ starting at 1 and the total time span fixed at 1. Applying the two models to the above examples, we estimate volatility processes.

Data are generated by the Euler method with data generating time interval $1/1.28\times 10^6$. On the other hand, observations are sampled out of them depending on sampling time interval $\Delta t=1/4,000$,
$1/8,000$, and $1/16,000$. We set the first $1/2$ period as the burn-in time in order to avoid the influence of the starting value of $X$. And then, the subsequent 1 period is used for estimation except that the first $1/40$ period is used for estimating the initial state of $\tilde Y$ by the sum of squared differences of $X$;
the initial state is given as $\sum_{k=1}^m(X_{s_k}-X_{s_{k-1}})^2/\Delta t$ for $\{X_{s_k}\}_{0\le k\le m}$, where $m$ depends on $\Delta t$ since the period for estimation of the initial state is fixed at $1/40$.

From discrete time series given as above, we estimate $\alpha$ and $\beta$ by LSE and $\theta$ by the quasi-MLE, and $\beta_0$ and $\beta_1$ by the least square estimation for the local linear model.
In the local linear model, the bandwidth $h$ is chosen by visual inspection.
We take $h=0.15$ for the first example, 0.13 for the second, 0.12 for the third, and 0.10 for the last.
Then, using their estimates, we compute $\tilde Y_{t|t}$ $(t\in\{t_k\}_{1\le k\le n})$ for the proposed model and $\hat m$ for the counterpart.

Figure 1 through 4 display the results of estimation of the volatility processes.
The left column shows the true curves of the volatility functions in a solid line and the estimated ones by the proposed model in a dotted line.
The right column shows the true curves and the estimated ones by the local linear model in the same way as the left column.
It can be easily seen that the estimated curves are converging to the true ones as $\Delta t$ becomes shorter. Particularly, the convergence is pronounced in the proposed model. Comparing the estimated curves by the two models, the proposed model produces more smooth curves while somewhat wiggly for the local linear model. This wiggly curves might imply too small bandwidth, but the results are almost the same or worse in another choice.
Anyway, we could confirm numerically the consistency of the proposed model that is proved in the previous section.

\subsection{Estimation in out-of-sample}
To evaluate the estimates given by the proposed model, we compare its performance with those by the local linear model in out-of-sample manner. To this end, we simulate 1,000 sample paths with $\Delta t=1/16,000$ while data generating time interval $1/3.2\times 10^5$. For every sample path, we use the first 2,000 data for the parameter estimation, and then, estimate the volatility states for the last 2,000 data. The estimation error is measured with the root mean squared errors (RMSE) based on the 2,000 states. Using the sample mean and standard deviation of 1,000 RMSEs, we compare the performance of the two models.

Here we consider frequently used interest rate models as follows:
\begin{eqnarray*}
dX_t&=&(0.184-0.2146X_t)dt+0.0783\sqrt{X_t}dB_t\hspace{1.6cm}{\bf(lin)}\\
dX_t&=&(0.0073-0.1409X_t)dt+0.2596X_tdB_t\hspace{1.7cm}{\bf(quad)}\\
dX_t&=&(0.0408-0.5921X_t)dt+1.2924X_t^{1.5}dB_t\hspace{1.4cm}{\bf(cube)}\\
dX_t&=&(0.0074-0.1180X_t)dt+0.0713X_t^{0.7296}dB_t\hspace{1cm}{\bf(nlin)}
\end{eqnarray*}
The parameters of the first and fourth examples are cited from  Fan and Zhang (2003), the second ones from Takamizawa and Shoji (2004), and the third ones from Chan et al (1992). The first example has a linear volatility, the second is quadratic, the third is cubic and the last is nonlinear.

Data are generated as the starting value $X_0=0.1$. But, like the previous experiment, the first 2,000 data are discarded in order to get rid of the influence of choice of the starting value. The results are presented in Table 1. Except for ({\bf lin}), the proposed model {\bf(semi)} show better performance in mean than the local linear model {\bf(ker)}. Particularly, looking at the standard deviations of RMSEs, the stable performance of the proposed model is pronounced.

\subsection{Estimation of integrated volatility}
It is interesting to construct the integrated volatility from spot volatilities estimated in the previous section and to compare those with the realized volatility, $R$, given as $\sum_{k=1}^n|X_{t_k}-X_{t_{k-1}}|^2$. Here we use an approximation of the integrated volatility given as,
\begin{eqnarray*}
\int_s^t\sigma^2(X_u)du&\approx&\sum_{k=1}^n\hat\sigma^2(X_{t_{k-1}})(t_k-t_{k-1})\\
&=&\sum_{k=1}^n\hat\sigma^2(X_{t_{k-1}})\Delta t
\end{eqnarray*}
where $\hat\sigma$ stands for the estimate of the diffusion coefficient and
$s=t_0<t_1<\cdots<t_n=t$ with $\Delta t=t_k-t_{k-1}$.
 Let $V_{semi}$ and $V_{ker}$ be the approximate integrated volatilities computed from spot volatilities estimated by the proposed model and the local linear model, respectively. In the same setting as the previous section, we compute these integrated volatilities in out-of-sample manner. That is, $V_{semi}$, $V_{ker}$ and $R$ are computed from the last 2,000 data. Then, we get the difference between $R$ and $V_{semi}$ and between $R$ and $V_{ker}$ as $R-V_{semi}$ and $R-V_{ker}$ for each sample path. The total differences are measure with the mean and standard deviation of differences for 1,000 sample paths. The results are presented in Table 2. The total differences are almost the same between the proposed model and the local linear model. And, the two models underestimate the integrated volatility as compare with the realized volatility. Furthermore, looking at the standard deviations, the difference between $V_{semi}$ and $V_{ker}$ is quite small as compared with the estimation of the spot volatility. Considering the stable estimation by ({\bf semi}),
this maybe implies that the realized volatility is volatile enough to cancel out the difference of the two models.

\section{Concluding remarks}
The paper proposed a semiparametric model of estimating the volatility defined by the squared diffusion coefficient of a stochastic differential equation. The volatility was approximated by a second order polynomial with stochastic coefficients and thereby we set up a vector process consisting of observable and unobservable processes in which the volatility process is defined as an unobservable one. By using the recursive updating formula, the volatility processes could be estimated by the filtering.

From theoretical viewpoints, we presented the proof of consistency of the proposed model in the sense that estimated processes converge to the true ones as the sampling interval goes to zero while the total time span goes to infinity.

And, from numerical viewpoints, we carried out the Monte Carlo experiments by which we could well estimate unobservable volatility processes and, at the same time, we confirmed the consistency numerically by using stochastic differential equations with linear/nonlinear diffusion coefficients.
Furthermore, through the performance comparison with the local linear model, the propose model showed better performance of volatility estimation in mean and standard deviation of estimation errors than the local linear model.


\setcounter{theorem}{1}
\section{Appendix for proofs}

In the following, $E_s[\cdot]$ stands for $E[\cdot|\F_s]$.
\begin{lemma}\label{lm:EXm}
For $s< t$ with $\Delta t=t-s$ and any positive integer $m$, $E_s|X_t-X_s|^{2m}$
denoted by $O((\Delta t)^m)$. That is, there is a constant $K_m$ depending on $m$ such that,
\begin{eqnarray*}
E_s|X_t-X_s|^{2m}&\le& K_m(\Delta t)^m,
\end{eqnarray*}
for sufficiently small $\Delta t$.
\end{lemma}
\noindent {\bf proof:}
$\mu(\cdot)$ stands for $\mu(\cdot;\eta)$ in this proof.
We prove the lemma by induction on $m$. Firstly, consider $m=1$. By the It\^o's 
formula,
\begin{eqnarray*}
(X_t-X_s)^2&=&
2\int_s^t\mu(X_u)(X_u-X_s)du
+\int_s^tg_udu\\
&&+2\int_s^t(X_u-X_s)\sqrt{g_u}dB_u\\
&\le&
\int_s^t(\mu(X_u)^2+(X_u-X_s)^2)du
+\int_s^tg_udu\\
&&+2\int_s^t(X_u-X_s)\sqrt{g_u}dB_u\\
&\le&
\int_s^tL_1du
+\int_s^t(X_u-X_s)^2du
+2\int_s^t(X_u-X_s)\sqrt{g_u}dB_u
\end{eqnarray*}
where $L_1$ stands for some constant since $X_u$ is bounded, 
and so are $\mu(X_u)$ and $g(X_u)$. Applying the conditional expectation at time $s$,
\[
E_s|X_t-X_s|^2\le
L_1\Delta t+\int_s^tE_s|X_u-X_s|^2du.
\]
By the Gronwall inequality, see Karatzas and Shreve (1991) for example, we get,
\[
E_s|X_t-X_s|^2\le
L_1\Delta t+\int_s^tL_1(u-s)e^{t-u}du.
\]
The second term has the order of $(\Delta t)^2$.
Actually, suppose
$\lim_{\Delta t\to 0}\int_s^t(u-s)e^{t-u}du/(\Delta t)^2$. Since,
\[
\int_s^t(u-s)e^{t-u}du
=
\int_0^{\Delta t}ue^{\Delta t-u}du
\]
we get,
\begin{eqnarray*}
\lim_{\Delta t\to 0}{\int_s^t(u-s)e^{t-u}du\over (\Delta t)^2}
&=&
\lim_{\Delta t\to 0}{\int_0^{\Delta t}ue^{\Delta t-u}du\over (\Delta 
t)^2}\\
&=&
\lim_{\Delta t\to 0}{\Delta t\over 2\Delta t}
\end{eqnarray*}
The claim holds for $m=1$.

Next, Suppose $E_s|X_t-X_s|^{2m}$. By the It\^o's formula,
\begin{eqnarray*}
(X_t-X_s)^{2m}&=&
2m\int_s^t\mu(X_u)(X_u-X_s)^{2m-1}du
+m(2m-1)\int_s^t(X_u-X_s)^{2(m-1)}g_udu\\
&&+2m\int_s^t(X_u-X_s)^{2m-1}\sqrt{g_u}dB_u\\
&\le&
m\int_s^t(\mu(X_u)^2+|X_u-X_s|^2)(X_u-X_s)^{2(m-1)}du\\
&&+m(2m-1)\int_s^t(X_u-X_s)^{2(m-1)}g_udu
+2m\int_s^t(X_u-X_s)^{2m-1}\sqrt{g_u}dB_u\\
&\le&
2m^2L_m\int_s^t(X_u-X_s)^{2(m-1)}du
+2m\int_s^t(X_u-X_s)^{2m}du\\
&&+2m\int_s^t(X_u-X_s)^{2m-1}\sqrt{g_u}dB_u
\end{eqnarray*}
for some constant $L_2$ such that $\mu(X_u)^2<L_2$ and $g_u<L_2$ since 
$X_s$ and $g_u$ are bounded. Hence,
\[
E_s|X_t-X_s|^{2m}
\le
2m^2L_2\int_s^tE_s|X_u-X_s|^{2(m-1)}du
+2m\int_s^tE_s|X_u-X_s|^{2m}du
\]
By the induction, there is a constant $K_{m-1}$ such that,
\begin{eqnarray*}
E_s|X_t-X_s|^{2m}
&\le&
2m^2L_2\int_s^tK_{m-1}(u-s)^{m-1}du
+2m\int_s^tE_s|X_u-X_s|^{2m}du\\
&=&
2m^2L_2K_{m-1}(\Delta t)^m
+2m\int_s^tE_s|X_u-X_s|^{2m}du
\end{eqnarray*}
By the Gronwall inequality,
\[
E_s|X_t-X_s|^{2m}
\le
2m^2L_2K_{m-1}(\Delta t)^m
+2m\int_s^t2m^2L_2K_{m-1}(u-s)^me^{2m(t-u)}du
\]
The claim holds for $m$. This completes the proof.

\begin{lemma}\label{lm:XaX}
For any $k$ $(1\le k\le n)$, let $t$ and $s$ be $t_k$ and $t_{k-1}$, respectively.
The order of $E_s|X_t-\tilde X_t|^2$ is $O(\Delta t)$.
\end{lemma}
\noindent {\bf proof:}
In this proof, $\mu(\cdot)$ and $\mu'(\cdot)$ stand for $\mu(\cdot;\eta)$ and $\mu'(\cdot;\eta)$, respectively.
By the It\^o's formula,
\begin{eqnarray*}
d(X-\tilde X)^2&=&2(X-\tilde X)dX-2(X-\tilde X)d\tilde 
X+d\<X\>+d\<\tilde X\>-2d\<X,\tilde X\>\\
&=&2(X-\tilde X)(dX-d\tilde X)+(\sqrt{g}-\sqrt{\tilde Y})^2dt\\
&=&2(X-\tilde X)(\mu(X;\eta)-\tilde\mu(\tilde X;\hat\eta))dt
+2(X-\tilde X)(\sqrt{g}-\sqrt{\tilde Y})dB
+(\sqrt{g}-\sqrt{\tilde Y})^2dt
\end{eqnarray*}
For simplicity, we may omit time subscription unless otherwise confusion.
Here, consider the first order Taylor expansion of $\mu$ and $\tilde\mu$. For $\mu$, there exist
$\nu\in[s,t]$ such that
$\mu(X_t)=\mu(X_s)+\mu'(X_\nu)(X_t-X_s)$.
For $\tilde\mu$, take the expansion around $\eta$. So,
$\tilde\mu(\tilde X_t;\hat\eta)=\tilde\mu(\tilde X_t;\eta)+\partial_\eta\tilde\mu(\tilde X_t;\bar\eta)(\hat\eta-\eta)$ for some $\bar\eta$, where $\partial_\eta$ stands for the gradient of $\tilde\mu$. The rate of convergence of $\hat\eta$ is $(N\Delta t)^{-{1\over 2}}$. Here note $\tau=n\Delta t$.
Hence,
\begin{eqnarray}\label{df:drft}
\mu(X;\eta)-\tilde\mu(\tilde X;\hat\eta)
&=&
\mu(X;\eta)-\tilde\mu(\tilde X;\eta)-\partial_\eta\tilde\mu(\tilde X;\bar\eta)(\hat\eta-\eta)\\
&=&
\mu(X_s)+\mu'(X_\nu)(X-X_s)-\{\mu(\tilde X_s)+\mu'(\tilde X_s)(\tilde X-\tilde X_s)\}\nonumber\\
&&-\partial_\eta\tilde\mu(\tilde X;\bar\eta)(\hat\eta-\eta)\nonumber\\
&=&\mu'(X_s)(X-\tilde X)+(\mu'(X_\nu)-\mu'(X_s))(X-X_s)
-\partial_\eta\tilde\mu(\tilde X;\bar\eta)(\hat\eta-\eta)\nonumber
\end{eqnarray}
Here note $X_s=\tilde X_s$ by the setting.
Using this,
\begin{eqnarray*}
d(X-\tilde X)^2
&=&
2\mu'(X_s)(X-\tilde X)^2dt
+2(X-\tilde X)(\mu'(X_\nu)-\mu'(X_s))(X-X_s)dt\\
&&-2(X-\tilde X)\partial_\eta\tilde\mu(\tilde X;\bar\eta)(\hat\eta-\eta)
+2(X-\tilde X)(\sqrt{g}-\sqrt{\tilde Y})dB
+(\sqrt{g}-\sqrt{\tilde Y})^2dt
\end{eqnarray*}
$\mu'<L$ because of the Lipschitz condition and both $g$ and $\tilde Y$ are all bounded. So, for some  $K_1$, $\mu'<K_1$, $g<K_1$ and $\tilde Y<K_1$. And note $\Delta t<\tau$. Then,
\begin{eqnarray*}
E_s|X_t-\tilde X_t|^2
&\le&
2K_1\int_s^tE_s|X_u-\tilde X_u|^2du
+4K_1\int_s^tE_s|(X_u-\tilde X_u)(X_u-X_s)|du\\
&&
+2\int_s^tE_s|(X_u-\tilde X_u)\partial_\eta\tilde\mu(\tilde X_u;\bar\eta)(\hat\eta-\eta)|du+2K_1\Delta t
\end{eqnarray*}
First, by using lemma \ref{lm:EXm},
\begin{eqnarray*}
2\int_s^tE_s|(X_u-\tilde X_u)(X_u-X_s)|du
&\le&
\int_s^tE_s|X_u-\tilde X_u|^2du
+\int_s^tE_s|X_u-X_s|^2du\\
&\le&
\int_s^tE_s|X_u-\tilde X_u|^2du
+K_2(\Delta t)^2,
\end{eqnarray*}
for some constant $K_2$.
Next, since $E_s|\partial_\eta\tilde\mu(\tilde X;\bar\eta)(\hat\eta-\eta)|^2<K_3(N\Delta t)^{-1}$ for some constant $K_3$,
\begin{eqnarray*}
2\int_s^tE_s|(X_u-\tilde X_u)\partial_\eta\tilde\mu(\tilde X_u;\bar\eta)(\hat\eta-\eta)|du
&\le&
\int_s^tE_s|X_u-\tilde X_u|^2du
+
\int_s^tE_s|\partial_\eta\tilde\mu(\tilde X_u;\bar\eta)(\hat\eta-\eta)|^2du\\
&\le&
\int_s^tE_s|X_u-\tilde X_u|^2du
+
K_3/N\\
&\le&
\int_s^tE_s|X_u-\tilde X_u|^2du
+
K_3\Delta t.
\end{eqnarray*}
We get the last inequality from $1/N<T/N=\Delta t$ since we consider $T\to\infty$.
Hence,
\begin{eqnarray*}
E_s|X_t-\tilde X_t|^2
&\le&
(4K_1+1)\int_s^tE_s|X_u-\tilde X_u|^2du
+2K_1\Delta t+2K_1K_2(\Delta t)^2+K_3\Delta t\\
&\le&
(4K_1+1)\int_s^tE_s|X_u-\tilde X_u|^2du
+(2K_1+2K_1K_2\tau+K_3)\Delta t.
\end{eqnarray*}
By the Gronwall inequality,
\[
E_s|X_t-\tilde X_t|^2\le
(2K_1+2K_1K_2\tau+K_3)\Delta t+(4K_1+1)\int_s^t(2K_1+2K_1K_2\tau+K_3)(u-s)e^{(4K_1+1)(t-u)}du.
\]
This completes the proof.

\begin{lemma}\label{lm:gY}
There is a positive constant $K$ such that,
\[
(g_t-Y_t)^2\le K((X_t-X_s)^2+(X_t-X_s)^4)
\]
where $t$ and $s$ stand for $t_k$ and $t_{k-1}$, respectively.
\end{lemma}
\noindent
{\bf proof:} The second order Taylor expansion of $g_t\equiv g(X_t)$ around $X_s$ is 
given as,
\[
g_t=g_s+g'(X_s)(X_t-X_s)+{g''(X_\eta)\over 2}(X_t-X_s)^2
\]
where $X_\eta=(1-\eta)X_s+\eta X_t$ for some $\eta\in [0,1]$.
From (\ref{def:f}) we have,
\begin{eqnarray*}
Y_t&=&Y_s+Y^1_s(X_t-X_s)+{\theta\over 2}(X_t-X_s)^2\\
Y_s&=&g(X_s)
\end{eqnarray*}
Since $X_u$ is a bounded process, $g'(X_u)$, $g''(X_u)$ and $Y_u$ are all bounded. 
Hence,
\begin{eqnarray*}
(g_t-Y_t)^2&=&
\left\{
(g'(X_s)-Y^1_s)(X_t-X_s)
+{1\over 2}(g''(X_\eta)-\theta)(X_t-X_s)^2
\right\}^2\\
&\le&
2\left\{
(g'(X_s)-Y^1_s)^2(X_t-X_s)^2
+\left({1\over 2}(g''(X_\eta)-\theta)\right)^2(X_t-X_s)^4
\right\}\\
&\le&
K(
(X_t-X_s)^2+(X_t-X_s)^4)
\end{eqnarray*}
for some positive constant $K$.


\medskip
\noindent {\bf proof of theorem \ref{th:YaY}:}
$\mu(\cdot)$ and $\mu'(\cdot)$ stand for $\mu(\cdot;\eta)$ and $\mu'(\cdot;\eta)$, respectively.
By the It\^o's formula,
\begin{eqnarray*}
d(Y-\tilde Y)^2
&=&
2(Y-\tilde Y)(dY-d\tilde Y)+d\<Y\>+d\<\tilde Y\>-2d\<Y,\tilde Y\>\\
&=&
2(Y-\tilde Y)
(Y^1dX-Y^1_sd\tilde X+{\theta\over 2}(g-\tilde Y)dt)
+(Y^1\sqrt{g}-Y^1_s\sqrt{\tilde Y})^2dt
\end{eqnarray*}
Here, we denote processes at time $s$ by $X_s$ for example. Firstly,
\[
Y^1dX-Y^1_sd\tilde X
=
(Y^1\mu(X)-Y^1_s\tilde\mu(\tilde X;\hat\eta))dt
+(Y^1\sqrt{g}-Y^1_s\sqrt{\tilde Y})dB
\]
Using $Y^1_t=Y^1_s+\theta(X_t-X_s)$ and (\ref{df:drft}), the coefficient of $dt$ is given as,
\begin{eqnarray*}
Y^1\mu(X)-Y^1_s\tilde\mu(\tilde X;\hat\eta)
&=&Y^1_s(\mu(X)-\tilde\mu(\tilde X;\hat\eta))+\theta\mu(X)(X-X_s)\\
&=&
Y^1_s\{
\mu'(X_s)(X-\tilde X)+
(\mu'(X_\nu)-\mu'(X_s))(X-X_s)
-\partial_\eta\tilde\mu(\tilde X;\bar\eta)(\hat\eta-\eta)\}\\
&&+\theta\mu(X)(X-X_s)\\
&=&
Y^1_s\mu'(X_s)(X-\tilde X)+
\{Y^1_s(\mu'(X_\nu)-\mu'(X_s))+\theta\mu(X)\}(X-X_s)\\
&&-Y^1_s\partial_\eta\tilde\mu(\tilde X;\bar\eta)(\hat\eta-\eta)
\end{eqnarray*}
Hence,
\begin{eqnarray*}
d(Y-\tilde Y)^2
&=&
2(Y-\tilde Y)
\{Y^1_s\mu'(X_s)(X-\tilde X)
+(Y^1_s(\mu'(X_\nu)-\mu'(X_s))+\theta\mu(X))(X-X_s)\\
&&\phantom{2(Y-\tilde Y)\{}
-Y^1_s\partial_\eta\tilde\mu(\tilde X;\bar\eta)(\hat\eta-\eta)+{\theta\over 2}(g-\tilde Y)
\}dt\\
&&+(Y^1\sqrt{g}-Y^1_s\sqrt{\tilde Y})^2dt\\
&&+2(Y-\tilde Y)(Y^1\sqrt{g}-Y^1_s\sqrt{\tilde Y})dB\\
&\le&
(Y-\tilde Y)^2dt\\
&&+
\{
Y^1_s\mu'(X_s)(X-\tilde X)
+(Y^1_s(\mu'(X_\nu)-\mu'(X_s))+\theta\mu(X))(X-X_s)\\
&&\phantom{+\{}
-Y^1_s\partial_\eta\tilde\mu(\tilde X;\bar\eta)(\hat\eta-\eta)+{\theta\over 2}(g-\tilde Y)\}^2dt\\
&&+(Y^1\sqrt{g}-Y^1_s\sqrt{\tilde Y})^2dt\\
&&+2(Y-\tilde Y)(Y^1\sqrt{g}-Y^1_s\sqrt{\tilde Y})dB\\
&\le&
(Y-\tilde Y)^2dt\\
&&
+4\{
(Y^1_s\mu'(X_s))^2(X-\tilde X)^2
+(Y^1_s(\mu'(X_\nu)-\mu'(X_s))+\theta\mu(X))^2(X-X_s)^2\\
&&\phantom{+3\{}
+
(Y^1_s\partial_\eta\tilde\mu(\tilde X;\bar\eta)(\hat\eta-\eta))^2
+{\theta^2\over 4}(g-\tilde Y)^2
\}dt\\
&&+(Y^1\sqrt{g}-Y^1_s\sqrt{\tilde Y})^2dt\\
&&+2(Y-\tilde Y)(Y^1\sqrt{g}-Y^1_s\sqrt{\tilde Y})dB
\end{eqnarray*}
$Y^1_s$, $\mu$, and $\mu'$ are all bounded. And, the rate of convergence of $\hat\eta$ is $(N\Delta t)^{-{1\over 2}}$. Hence, without loss of generality, there is a constant $K_1$ such that,
\begin{eqnarray*}
d(Y-\tilde Y)^2
&\le&
(Y-\tilde Y)^2dt
+
K_1((X-\tilde X)^2+(X-X_s)^2+(g-\tilde Y)^2+1/(N\Delta t))dt\\
&&+(Y^1\sqrt{g}-Y^1_s\sqrt{\tilde Y})^2dt\\
&&+2(Y-\tilde Y)(Y^1\sqrt{g}-Y^1_s\sqrt{\tilde Y})dB
\end{eqnarray*}

Firstly,
\begin{eqnarray*}
(g_t-\tilde Y_t)^2&=&\{(g_t-Y_t)+(Y_t-\tilde Y_t)\}^2\\
&\le&
2((g_t-Y_t)^2+(Y_t-\tilde Y_t)^2)
\end{eqnarray*}
Thanks to lemma \ref{lm:gY}, there is a positive constant $K_2$ such that,
\[
(g_t-Y_t)^2\le K_2((X_t-X_s)^2+(X_t-X_s)^4)
\]

Next, we want to evaluate the coefficient of $dt$ in the second line.
\begin{eqnarray*}
(Y^1\sqrt{g_t}-Y^1_s\sqrt{\tilde Y_t})^2
&=&
(Y^1_s(\sqrt{g_t}-\sqrt{\tilde Y_t})+\theta\sqrt{g_t}(X_t-X_s))^2\\
&\le&
2\{(Y^1_s(\sqrt{g_t}-\sqrt{\tilde Y_t}))^2
+(\theta\sqrt{g_t}(X_t-X_s))^2\}
\end{eqnarray*}
Since $X$ is a bounded process, $Y^1$ and $g$ are also bounded. Hence, there is a 
positive constant $K_3$ such that,
\[
(Y^1\sqrt{g_t}-Y^1_s\sqrt{\tilde Y_t})^2
\le
K_3((\sqrt{g_t}-\sqrt{\tilde Y_t})^2
+(X_t-X_s)^2)
\]

Furthermore
\begin{eqnarray*}
(\sqrt{g_t}-\sqrt{\tilde Y_t})^2
&=&
\{
(\sqrt{g_t}-\sqrt{g_s})
+(\sqrt{g_s}-\sqrt{Y_t})
+(\sqrt{Y_t}-\sqrt{\tilde Y_t})
\}^2\\
&\le&
3\{
(\sqrt{g_t}-\sqrt{g_s})^2
+(\sqrt{g_s}-\sqrt{Y_t})^2
+(\sqrt{Y_t}-\sqrt{\tilde Y_t})^2
\}
\end{eqnarray*}
Because of the Lipschitz condition of $\sqrt{g}$,
\[
(\sqrt{g_t}-\sqrt{g_s})^2\le L^2(X_t-X_s)^2
\]

To evaluate the second and third terms, we introduce stopping times for a sufficiently 
small $\epsilon>0$ as follows:
\begin{eqnarray*}
T^\epsilon_3&=& \inf\{t\ge 0; \tilde Y_t\le \epsilon\}.
\end{eqnarray*}
As for $Y_t$, since its sample path isn't necessarily continuous, we firstly define a stopping time $S_k$ for $t\in[t_{k-1},t_k)$ as follows:
\[
S_k=\inf\{t\ge t_{k-1}; Y_t\le \epsilon\}.
\]
And then, a stopping time $S$ is defined by,
\[
S=
\left\{
\begin{tabular}{ll}
$S_k$& $t\in[t_{k-1},t_k)$ $(1\le k\le n)$\\
$\tau$& $t\ge t_n=\tau$
\end{tabular}
\right.
\]
Using these stopping times, we newly redefine $Y_t$ and $\tilde Y_t$ as
$Y_t\equiv Y_{t\katsu T^\epsilon_3\katsu S}$ and
$\tilde Y_t\equiv \tilde Y_{t\katsu T^\epsilon_3\katsu S}$, 
respectively.

\begin{eqnarray*}
(\sqrt{g_s}-\sqrt{Y_t})^2
&=&
(\sqrt{Y_t}-\sqrt{Y_s})^2\\
&=&
\left(
{Y_t-Y_s\over \sqrt{Y_t}+\sqrt{Y_s}}
\right)^2\\
&\le&
{(Y_t-Y_s)^2\over Y_t+Y_s}\\
&\le&
{(Y_t-Y_s)^2\over 2\epsilon}
\end{eqnarray*}

Similarly,
\begin{eqnarray*}
(\sqrt{Y_t}-\sqrt{\tilde Y_t})^2
&=&
\left(
{Y_t-\tilde Y_t\over \sqrt{Y_t}+\sqrt{\tilde Y_t}}
\right)^2\\
&\le&
{(Y_t-\tilde Y_t)^2\over Y_t+\tilde Y_t}\\
&\le&
{(Y_t-\tilde Y_t)^2\over 2\epsilon}
\end{eqnarray*}

Using the above inequalities, we get,
\begin{eqnarray*}
d(Y-\tilde Y)^2
&\le&
K_1(X-\tilde X)^2dt+K_4(X-X_s)^2dt+K_5(X-X_s)^4dt+K_6(Y-\tilde Y)^2dt\\
&&+K_1/(N\Delta t)dt+2(Y-\tilde Y)(Y^1\sqrt{g}-Y^1_s\sqrt{\tilde Y})dB
\end{eqnarray*}
where,
\begin{eqnarray*}
K_4&=&K_1(1+2K_2)+K_3(1+3L^2+{3K_2\over 2\epsilon})\\
K_5&=&2K_1K_2+{3K_2K_3\over 2\epsilon}\\
K_6&=&1+2K_1+{3K_3\over 2\epsilon}
\end{eqnarray*}
Hence,
\begin{eqnarray*}
E_s|Y_t-\tilde Y_t|^2-E_s|Y_s-\tilde Y_s|^2
&=&
E_s|Y_t-\tilde Y_t|^2-|Y_s-\tilde Y_s|^2\\
&\le&
K_1\int_s^tE_s|X_u-\tilde X_u|^2du
+
K_4\int_s^tE_s|X_u-X_s|^2du\\
&&+
K_5\int_s^tE_s|X_u-X_s|^4du
+
K_6\int_s^tE_s|Y_u-\tilde Y_u|^2du+1/N.
\end{eqnarray*}
Thanks to lemma \ref{lm:EXm} and \ref{lm:XaX}, the integrands of the first and 
second terms have the order of $(u-s)$ and the third $(u-s)^2$. Consequently their integrations have the order of $(\Delta t)^2$ and $(\Delta t)^3$.
The resulting inequality is give by,
\[
E_s|Y_t-\tilde Y_t|^2
\le
|Y_s-\tilde Y_s|^2+\phi(t)+K_6\int_s^tE_s|Y_u-\tilde Y_u|^2du
\]
where $\phi(t)$ consists of the two parts; one has the order of $1/N$ and the othre has the order of $(\Delta t)^2$. By the Gronwall inequality, we get,
\begin{eqnarray*}
E_s|Y_t-\tilde Y_t|^2
&\le&
|Y_s-\tilde Y_s|^2
+\phi(t)+K_6\int_s^t(|Y_s-\tilde Y_s|^2
+\phi(u))e^{K_6(t-u)}du\\
&=&
|Y_s-\tilde Y_s|^2e^{K_6\Delta t}+\phi(t)+K_6\int_s^t\phi(u)e^{K_6(t-u)}du
\end{eqnarray*}
Note that the summation of the second and third terms, denote by $\psi$, can be expressed by $\phi_1(1/N)+\phi_2((\Delta t)^2)+\phi_3(\Delta t/N,(\Delta t)^3)$, where $\phi_1$ has the order of
$1/N$, $\phi_2$ has the $(\Delta t)^2$ order, and $\phi_3$ has the sum of $\Delta t/N$ and $(\Delta t)^3$.
Applying the unconditional expectation to the both sides and substituting $t$ and $s$ for $t_k$ and $t_{k-1}$, we get,
\[
E|Y_{t_k}-\tilde Y_{t_k}|^2
\le
e^{K_6\Delta t}E|Y_{t_{k-1}}-\tilde Y_{t_{k-1}}|^2+\psi
\]
Multiplying the both sides by $e^{(n-k)K_6\Delta t}$ and summing it up from $k=1$ to 
$n$,
\[
E|Y_{t_n}-\tilde Y_{t_n}|^2
\le
e^{(n-1)K_6\Delta t}E|Y_0-\tilde Y_0|^2
+{\psi(e^{nK_6\Delta t}-1)\over e^{K_6\Delta t}-1}
\]
But, $Y_0=\tilde Y_0$ by the setting. Hence
\[
E|Y_\tau-\tilde Y_\tau|^2
\le
{\psi(e^{nK_6\Delta t}-1)\over e^{K_6\Delta t}-1}
=
{\psi\over\Delta t}\cdot \Delta t{(e^{nK_6\Delta t}-1)\over e^{K_6\Delta t}-1}
\]
On one hand, since $n\Delta t=\tau$,
\begin{eqnarray*}
\lim_{\Delta t\to 0}\Delta t{(e^{nK_6\Delta t}-1)\over e^{K_6\Delta t}-1}
&=&
\lim_{\Delta t\to 0}{(e^{K_6\tau}-1)\Delta t\over e^{K_6\Delta t}-1}\\
&=&
{e^{K_6\tau}-1\over K_6}
\end{eqnarray*}
On the other hand, $\psi/\Delta t$ can be expressed by
$\phi_1(1/(N\Delta t))+\phi_2(\Delta t)+\phi_3(1/N,(\Delta t)^2)$.
Hence $E|Y_\tau-\tilde Y_\tau|^2$ converges to zero as $N\to\infty$, $\Delta t\to0$ and $N\Delta t\to\infty$.
Finally, by letting
$\epsilon\downarrow 0$, we can get the desired result.

\begin{lemma}\label{lm:V12}
Let $V_1(x)$ and $V_2(x)$ be the (1,1) and (1,2) elements of $V_{t|s}(x)$. For any 
$\epsilon>0$, there exists a $\delta>0$ such that for any $\Delta t<\delta$,
\[
V_1(x)\left|
{V_2(x)\over V_1(x)}-{V_2(y)\over V_1(y)}
\right|^2
\le
\epsilon\Delta t
\]
\end{lemma}
\noindent
{\bf proof:} Now suppose $V_i(x)$ $(i=1,2)$ is a function of $\Delta t$, which is 
denoted by $v_i(\Delta t;x)$.
Then, we want to evaluate
\[
\lim_{\Delta t\to0}
{v_2(\Delta t;x)\over v_1(\Delta t;x)}-{v_2(\Delta t;y)\over v_1(\Delta t;y)}
\]
Here note $\lim_{\Delta t\to0}v_i(\Delta t;x)=0$
and $\lim_{\Delta t\to0}v_i'(\Delta t;x)=x$.
For simplicity, we denote $v_i(\Delta t;x)$ and $v_i(\Delta t;y)$ by
$x_i$ and $y_i$, respectively.
\begin{eqnarray*}
\lim_{\Delta t\to0}
{x_2\over x_1}-{y_2\over y_1}
&=&
\lim_{\Delta t\to0}
{x_2y_1-x_1y_2\over x_1y_1}\\
&=&
\lim_{\Delta t\to0}
{(x_2''y_1+2x_2'y_1'+x_2y_1'')-(x_1''y_2+2x_1'y_2'+x_1y_2'')
\over x_1''y_1+2x_1'y_1'+x_1y_1''}\\
&=&0
\end{eqnarray*}
On the other hand,
\begin{eqnarray*}
\lim_{\Delta t\to0}
{x_1\over\Delta t}
&=&
\lim_{\Delta t\to0}
{x_1'\over 1}\\
&=&x
\end{eqnarray*}
Hence,
\[
\lim_{\Delta t\to0}
{x_1\over\Delta t}
\left|
{x_2\over x_1}-{y_2\over y_1}
\right|^2
=0
\]
In other words, for any $\epsilon>0$, there exists a $\delta>0$ such that, for all 
$\Delta t<\delta$,
\[
{x_1\over\Delta t}
\left|
{x_2\over x_1}-{y_2\over y_1}
\right|^2
<\epsilon
\]
This completes the proof.

\begin{lemma}\label{lm:V123}

Let $V_1(x)$ and $V_2(x)$ be as above. And let $V_3(x)$ be the (2,2) elements of $V_{t|s}(x)$.
For any $\epsilon>0$, there exists a $\delta>0$ such that for all $\Delta t<\delta$,
\[
\left|
{V_1(x)V_3(x)-V_2(x)^2\over V_1(x)}
\right|
\le
\epsilon\Delta t
\]
\end{lemma}
\noindent
{\bf proof:} Now suppose $V_i(x)$ $(1\le i\le 3)$ are function of $\Delta t$, which 
is simply denoted by $x_i$. And note $\lim_{\Delta t\to0}x_i=0$,
$\lim_{\Delta t\to0}x_i'=x$ and $\lim_{\Delta t\to0}|x_i''|<\infty$.
\begin{eqnarray*}
\lim_{\Delta t\to0}
{x_1x_3-x_2^2\over x_1\Delta t}
&=&
\lim_{\Delta t\to0}
{(x_1''x_3+2x_1'x_3'+x_1x_3'')-2((x_2')^2+x_2x_2'')\over x_1''\Delta t+2x_1'}\\
&=&0
\end{eqnarray*}
Hence,
for any $\epsilon>0$, there exists a $\delta>0$ such that, for all $\Delta t<\delta$,
\[
\left|
{x_1x_3-x_2^2\over x_1\Delta t}
\right|
<\epsilon
\]
This completes the proof.

\noindent
\medskip
{\bf proof of theorem \ref{th:aYfY}:} For simplicity, let $t$ and $s$ be $t_k$ and $t_{k-1}$, respectively.
And, let $V_1(x)$, $V_2(x)$ and $V_3(x)$ be the (1,1), (1,2) and (2,2) elements of
$V_{t|s}(x)$.
Recall,
\[
\tilde Y_{t|t}=\tilde Y_{t|s}+\kappa(\tilde X_t-\tilde X_{t|s}),
\]
where,
\[
\kappa={V_2(\tilde Y_{s|s})\over V_1(\tilde Y_{s|s})}
\]
Hence,
\begin{eqnarray*}
E_s[(\tilde Y_t-\tilde Y_{t|t})^2]
&=&
E_s[(\tilde Y_t-\tilde Y_{t|s}-\kappa(\tilde X_t-\tilde X_{t|s}))^2]\\
&=&
E_s[(\tilde Y_t-\tilde Y_{t|s})^2]
-2\kappa
E_s[(\tilde Y_t-\tilde Y_{t|s})(\tilde X_t-\tilde X_{t|s})]^2
+\kappa^2
E_s[(\tilde X_t-\tilde X_{t|s})^2]
\end{eqnarray*}
Here note $\kappa\in\G_s$.

Firstly, we evaluate the first term. Noticing $\G_s\subset\F_s$,
\begin{eqnarray*}
E_s[(\tilde Y_t-\tilde Y_{t|s})^2]
&=&
E_s[\{(\tilde Y_t-E_s[\tilde Y_t])+(E_s[\tilde Y_t]-\tilde Y_{t|s})\}^2]\\
&=&
E_s[(\tilde Y_t-E_s[\tilde Y_t])^2]
+
E_s[(\tilde Y_t-E_s[\tilde Y_t])E_s[\tilde Y_t-\tilde Y_{t|s}]]
+
E_s[(E_s[\tilde Y_t]-\tilde Y_{t|s})^2]\\
&=&
V_3(\tilde Y_s)
+
e^{2c\Delta t}(\tilde Y_s-\tilde Y_{s|s})^2
\end{eqnarray*}

Secondly, since $\tilde X_{t|s}=E_s[\tilde X_t]$,
\begin{eqnarray*}
E_s[(\tilde Y_t-\tilde Y_{t|s})(\tilde X_t-\tilde X_{t|s})]
&=&
E_s[(\tilde Y_t-E_s[\tilde Y_t])(\tilde X_t-E_s[\tilde X_t])]\\
&&+
E_s[E_s[\tilde Y_t-\tilde Y_{t|s}](\tilde X_t-E_s[\tilde X_t])]\\
&=&
V_2(\tilde Y_s)
\end{eqnarray*}

Hence,
\begin{eqnarray*}
E_s[(\tilde Y_t-\tilde Y_{t|t})^2]
&=&
V_3(\tilde Y_s)
+
e^{2c\Delta t}(\tilde Y_s-\tilde Y_{s|s})^2
-2
{V_2(\tilde Y_{s|s})\over V_1(\tilde Y_{s|s})}V_2(\tilde Y_s)
+
\left({V_2(\tilde Y_{s|s})\over V_1(\tilde Y_{s|s})}\right)^2
V_1(\tilde Y_s)\\
&=&
e^{2c\Delta t}(\tilde Y_s-\tilde Y_{s|s})^2
+
V_1(\tilde Y_s)
\left(
{V_2(\tilde Y_s)\over V_1(\tilde Y_s)}
-{V_2(\tilde Y_{s|s})\over V_1(\tilde Y_{s|s})}
\right)^2\\
&&+
{V_1(\tilde Y_s)V_3(\tilde Y_s)-V_2(\tilde Y_s)^2\over V_1(\tilde Y_s)}
\end{eqnarray*}
Due to lemma \ref{lm:V12} and \ref{lm:V123}, for any $\epsilon>0$, there exits 
a $\delta>0$ such that for all $\Delta t<\delta$,
\begin{eqnarray*}
V_1(\tilde Y_s)
\left(
{V_2(\tilde Y_s)\over V_1(\tilde Y_s)}
-{V_2(\tilde Y_{s|s})\over V_1(\tilde Y_{s|s})}
\right)^2
&<&\epsilon\Delta t\\
\left|
{V_1(\tilde Y_s)V_3(\tilde Y_s)-V_2(\tilde Y_s)^2\over V_1(\tilde Y_s)}
\right|
&<&\epsilon\Delta t
\end{eqnarray*}
Hence,
\[
E_s[(\tilde Y_t-\tilde Y_{t|t})^2]
<
e^{2c\Delta t}(\tilde Y_s-\tilde Y_{s|s})^2
+2\epsilon\Delta t
\]
Applying the unconditional expectation, we get,
\[
E[(\tilde Y_t-\tilde Y_{t|t})^2]
<
e^{2c\Delta t}E[(\tilde Y_s-\tilde Y_{s|s})^2]
+2\epsilon\Delta t
\]
Recall $t$ and $s$ stand for $t_k$ and $t_{k-1}$, respectively.
By multiplying $e^{2c\Delta t(n-k)}$ by the both side and summing it up from $k=1$ to 
$n$, we get,
\begin{eqnarray*}
E[(\tilde Y_{t_n}-\tilde Y_{t_n|t_n})^2]
&=&
E[(\tilde Y_\tau-\tilde Y_{\tau|\tau})^2]\\
&<&
e^{2c\Delta t(n-1)}E[(\tilde Y_0-\tilde Y_{0|0})^2]
+2\epsilon\Delta t
{e^{2cn\Delta t}-1\over e^{2c\Delta t}-1}\\
&=&
2\epsilon\Delta t
{e^{2c\tau}-1\over e^{2c\Delta t}-1}
\end{eqnarray*}
Here note $\tilde Y_0=\tilde Y_{0|0}$ by the setting, $t_n=\tau$ and
$n\Delta t=\tau$. By $\Delta t$ going to zero,
\[
\lim_{\Delta t\to 0}E[(\tilde Y_\tau-\tilde Y_{\tau|\tau})^2]
\le
\epsilon{e^{2c\tau}-1\over c}
\]
Since $\epsilon$ is arbitrarily given,
this completes the proof.

\section*{References}
\begin{description}
\item{}
A\"{\i}t-Sahalia, Y. (1996).
Nonparametric pricing of interest rate derivative securities.
{\em Econometrica} {\bf 64}, 527-560.

\item{}
Andersen, T. G., Bollerslev, T., Diebold, F. X. and Labys, P. (2003).
Modeling and forecasting realized volatility.
{\em Econometrica} {\bf 71}, 579-625.

\item{}
Andersen, T. G., Bollerslev, T. and Meddahi, N. (2004).
Analytical evaluation of volatility forecasts.
{\em International Economic Review} {\bf 45}, 1079-1110.

\item{}
Andersen, T. G., Bollerslev, T. and Meddahi, N. (2005).
Correcting the errors: volatility forecast evaluation using high-frequency data and 
realized  volatilities.
{\em Econometrica} {\bf 73}, 279-296.

\item{}
Anderson, B. D. O. and Moore, J. B. (1979).
{\em Optimal Filtering}.
Prentice-Hall: New Jersey.

\item{}
Bandi, F. M. and Phillips, P. C. B. (2003)
Fully nonparametric estimation of scalar diffusion models.
{\em Econometrica} {\bf 71}, 241-283.

\item{}
Barndorff-Nielsen, O. E. and Shephard, N. (2002).
Econometric analysis of realized volatility and its use in estimating stochastic 
volatility models.
{\em J. R. Statist. Soc.} B {\bf 64}, 253-280.

\item{}
Bali, T. G. and Wu, L. (2006).
A comprehensive analysis of the short-term
interest-rate dynamics.
{\em Jouranl of Banking \& Finance} {\bf 30}, 1269-1290.

\item{}
Barndorff-Nielsen, O. E. and Shephard, N. (2004).
Econometric analysis of realized covariation: High frequency based covariance, 
regression, and correlation in financial economics.
{\em Econometrica} {\bf 72}, 885-925.

\item{}
Campbell, J.Y., Lo, A. W. and MacKinlay, A. C. (1997).
{\em The Econometrics of Financial Markets}.
Princeton University Press: Princeton, New Jersey.

\item{}
Chan, K. C., Karolyi, G. A., Longstaff, F. A. and Sanders, A. B. (1992).
An empirical comparison of alternative models of the short-term interest rate.
{\em Journal of Finance} {\bf 47}, 1209-1227.

\item{}
Chapman, D. A. and Pearson, N. D. (2000).
Is the short rate drift actually nonlinear?
{\em Journal of Finance} {\bf 55}, 355-388.

\item{}
Deo, R., Hurvich, C. and Lu, Y. (2006).
Forecasting realized volatility using a long-memory stochastic volatility model: 
estimation, prediction and seasonal adjustment.
{\em Journal of Econometrics} {\bf 131}, 29-58.

\item{}
Engle, R. F. and Gallo, G. M.. (2006).
A multiple indicators model for volatility using intra-daily data.
{\em Journal of Econometrics} {\bf 131}, 3-27.

\item{}
Fan, J. and Gijbels, I. (1996).
{\em Local Polynomial Modelling and Its Applications}.
Chapman \& Hall: London.

\item{}
Fan, J. and Zhang,, C. (2003).
A reexamination of diffusion estimators with applications to financial model validation.
{\em Journal of the American Statistical Association} {\bf 98}, 118-134.

\item{}
Fan, J. and Yao, Q. (1998).
Efficient estimation of conditional variance functions in stochastic regression.
{\em Biometrika} {\bf 85}, 645-660.

\item{}
Florens-Zmirou, D. (1989).
Approximate discrete-time schemes for statistics of diffusion processes.
{\em Statistics} {\bf 20}, 547-557.

\item{}
Florens-Zmirou, D. (1993).
On estimating the diffusion coefficient from discrete observations.
{\em Journal of Applied Probability} {\bf 30}, 790-804.

\item{}
Ghysels, E., Santa-Clara, P. and Valkanov, R. (2006).
Predicting volatility: getting the most out of return data sampled at different 
frequencies.
{\em Journal of Econometrics} {\bf 131}, 59-95.

\item{}
Karatzas, I. and Shreve, S. E. (1991).
{\em Brownian Motion and Stochastic Calculus: 2nd ed}.
Springer: New York.

\item{}
Kessler, M. (1997).
Estimation of an ergodic diffusion from discrete observations.
{\em Scandinavian Journal of Statistics} {\bf 24}, 211-224.

\item{}
Prakasa Rao, B. L. S. (1983).
Asymptotic theory for non-linear least squares estimator for diffusion processes.
{\em Math. Operationsforsch. Statist. Ser. Stat.} {\bf 14}, 195-209.

\item{}
Stanton, R. (1997).
A nonparametric model of term structure dynamics and the market price of interest rate
risk.
{\em Journal of Finance} {\bf 52}, 1973-2002.

\item{}
Sun, L. (2003).
Nonlinear drift and stochastic volatility: An empirical investigation of short-term interest rate models.
{\em Journal of Financial Research} {\bf 26}, 389-404.

\item{}
Takamizawa, H. and Shoji, I. (2004).
On the accuracy of the local linear approximation for the term structure of interest rates.
{\em Quantitative Finance} {\bf 4}, 151-157.

\item{}
Thomakos, D. D., Wang, T. and Wille, L. T. (2002).
Modeling daily realized futures volatility with singular spectrum analysis.
{\em Physica A} {\bf 312}, 505-519.

\item{}
Yoshida, N. (1992).
Estimation for diffusion processes from discrete observation.
{\em Journal of Multivariate Analysis} {\bf 41}, 220-242.

\end{description}

\clearpage
\begin{table}
\begin{center}
\begin{tabular}{lc|cccc}
\hline
&& {\bf lin}& {\bf quad}& {\bf cube}& {\bf nlin}\\
\hline
{\bf semi}&
mean&  1.0083 & 0.9454 & 2.3153 & 0.2597 \\
&std&  1.0397 & 0.9621 & 2.6087 & 0.2544 \\
\hline
{\bf ker}&
mean&  0.9747 & 1.0691 & 2.6172 & 0.2828 \\
&std&  1.1800 & 1.3238 & 3.7640 & 0.3272 \\
\hline
\end{tabular}
\end{center}
\caption{Means (mean) and standard deviations (std) of 1,000 RMSE's of the proposed model {\bf(semi)} and the local linear model {\bf(ker)} are presented. Actual values should be multiplied by $10^{-4}$.}
\end{table}

\clearpage
\begin{table}
\begin{center}
\begin{tabular}{lc|cccc}
\hline
&& {\bf lin}& {\bf quad}& {\bf cube}& {\bf nlin}\\
\hline
$R-V_{semi}$&
mean&  5.3036 & 5.7694 & 14.3737 & 1.5188 \\
&std&  0.6111 & 1.3237 & 8.3288 & 0.1456 \\

\hline
$R-V_{ker}$&
mean&  5.3039 & 5.7744 & 14.3431 & 1.5208 \\
&std&  0.6096 & 1.3240 & 8.2681 & 0.1452 \\
\hline
\end{tabular}
\end{center}
\caption{Means (mean) and standard deviations (std) of the differences for 1,000 sample paths are presented. Actual values should be multiplied by $10^{-4}$.}
\end{table}

\end{document}